\documentclass[12pt]{amsart}

\usepackage{amssymb, amscd}

\newtheorem{Lemma}             {\bf Lemma}

\newtheorem{Theorem}    [Lemma]{\bf Theorem}

\newcommand{\ov}[1]{\overline{#1}}

\begin{document}
\title{On a result of Kiyota, Okuyama and Wada}
\author{ John Murray}
\address{Department of Mathematics,
         National University of Ireland,\newline
	 Maynooth, Co.~Kildare, Ireland}
\email{  John.Murray@maths.nuim.ie}
\date{   June 25, 2012}

\begin{abstract}
M. Kiyota, T. Okuyama and T. Wada recently proved that each $2$-block of a symmetric group $\Sigma_n$ contains a unique irreducible Brauer character of height $0$. We present a more conceptual proof of this result.
\end{abstract}
\maketitle

%%%%%%%%%%%%%%%%%%%%%%%%%%%%%%%%%%%%%%%%%%%%%%%%%%%%%%%%%%%%%%%%%
\section{Background on bilinear forms}\label{S:Fong}

According to the main result in \cite{KOW}, every $2$-block of the symmetric group $\Sigma_n$ has a unique irreducible Brauer character of height $0$. This is not true for an arbitrary $2$-block of a finite group. For example, let $B$ be a real non-principal $2$-block which is Morita equivalent to the group algebra of $A_4$ and which has a Klein-four defect group and a dihedral {\em extended defect group} (in the sense of \cite{G}). Then one can show that $B$ has three real irreducible Brauer characters of height $0$. The non-principal $2$-block of $((C_2\times C_2):C_9):C_2$ is of this type.

In this note we place the results of \cite{KOW} in a more general context using the approach to bilinear forms developed by R. Gow and W. Willems \cite{GW93}. We use results and notation from \cite{NT} for representation theory, from \cite{J} for symmetric groups, and from \cite{W} for bilinear forms in characteristic $2$.

Let $G$ be a finite group and let $(K,R,F)$ be a {\em $2$-modular system} for $G$. So $R$ is a complete discrete valuation ring with field of fractions $K$ of characteristic $0$, and residue field $R/J=F$ of characteristic $2$. Assume that $K$ contains a primitive $|G|$-th root of unity, and that $F$ is perfect. Then $K$ and $F$ are splitting fields for each subgroup of $G$.

The anti-isomorphism $g\mapsto g^{-1}$ on $G$ extends to an involutary $F$-algebra anti-automorphism $\sigma:FG\rightarrow FG$ called the {\em contragredient map}. Let $V$ be a right $FG$-module. The linear dual $V^*:={\rm Hom}_F(V,F)$ is considered as a right $FG$-module via $(f.x)(v):=f(v x^\sigma)$, for $f\in V^*$, $x\in FG$ and all $v\in V$. The Frobenius automorphism $\lambda\mapsto\lambda^2$ of the field $F$ induces an automorphism $(a_{ij})\mapsto(a_{ij}^2)$ of the group ${\rm GL}_F(V)$. Composing the module map $G\rightarrow{\rm GL}_F(V)$ with this automorphism endows $V$ with another $FG$-module structure. This module is called the {\em Frobenius twist} of $V$, and is denoted $V^{(2)}$.

Let $V^*\otimes V^*$ be the space of bilinear forms on $V$ and let $\Lambda^2(V^*)$ be the subspace of {\em symplectic bilinear forms} on $V$; a bilinear form $b:V\times V\rightarrow F$ is symplectic if and only if $b(v,v)=0$, for all $v\in V$. The quotient space $V^*\otimes V^*/\Lambda^2(V^*)$ is called the {\em symmetric square} of $V^*$ and is denoted $S^2(V^*)$.

A {\em quadratic form} on $V$ is a map $Q:V\rightarrow F$ such that $Q(\lambda v)=\lambda^2Q(v)$ and $(u,v)\mapsto Q(u+v)-Q(u)-Q(v)$ is a bilinear form on $V$, for all $u,v\in V$ and $\lambda\in F$. Now if $b$ is a bilinear form, its {\em diagonal} $\delta(b):v\mapsto b(v,v)$ is a quadratic form. The assignment $\delta$ is linear with kernel $\Lambda^2(V^*)$. So there is a short exact sequence of vector spaces
\begin{equation}\label{E:ses1}
\begin{CD}
0 @>>> \Lambda^2(V^*) @>>> V^*\otimes V^* @>\delta>> S^2(V^*) @>>> 0
\end{CD}
\end{equation}
We may identify $S^2(V^*)$ with the space of quadratic forms on $V$. If $Q$ is a quadratic form, its {\em polarization} is the associated bilinear form $\rho(Q):(u,v)\mapsto Q(u+v)-Q(u)-Q(v)$.

The dual $S^2(V)^*$ of the symmetric square $S^2(V)$ of $V$ is the space of {\em symmetric bilinear forms} on $V$. As char$(F)=2$, each symplectic form is symmetric. If $b$ is a symmetric bilinear form, $\delta(b)$ is additive and hence can be identified with a linear map $V^{(2)}\rightarrow F$. Thus there is a short exact sequence:
\begin{equation}\label{E:ses2}
\begin{CD}
0 @>>> \Lambda^2(V^*) @>>> S^2(V)^* @>\delta>> V^{(2)*} @>>> 0
\end{CD}
\end{equation}

All of these $F$-spaces are $FG$-modules, and the maps are $FG$-module homomorphisms. It is a singular feature of the characteristic $2$-theory that $S^2(V^*)$ and $S^2(V)^*$ need not be isomorphic as $FG$-modules.

Now let $b$ be a bilinear form on $V$. We say that $b$ is {\em $G$-invariant} if the associated map $v\mapsto b(v,\,\,)$ for $v\in V$, is an $FG$-module map $V\rightarrow V^*$. We say that $b$ is {\em nondegenerate} if this map is an $F$-isomorphism. Taking $G$-fixed points in \eqref{E:ses2} we get a long exact sequence of the form
$$
{\tiny
\begin{CD}
0\! @>>>\! \Lambda^2(V^*)^G\! @>>>\! S^2(V)^{*G}\!	@>\delta>>\! V^{(2)*G}\!	@>>>\! H^1(G,\Lambda^2(V^*))\! @>>>\dots\\
\end{CD}}
$$
In particular, if $V^{(2)*G}=0$, then each $G$-invariant symmetric bilinear form on $V$ is symplectic. Now the trivial $FG$-module equals its Frobenius twist. A simple argument then shows:

\begin{Lemma}\label{L:trivial}
If $V\cong V^*$, and $V$ has no trivial $G$-submodules, then each $G$-invariant symmetric bilinear form on $V$ is symplectic.
\end{Lemma}

We will make use of Fong's Lemma:

\begin{Lemma}
Let $V$ be an absolutely irreducible non-trivial $FG$-module. Then $V\cong V^*$ if and only if $V$ affords a nondegenerate $G$-invariant symplectic bilinear form. In particular dim$(V)$ is even.
\end{Lemma}

Let $t,h\in G$, with $t$ an involution and $h$ not an involution. Define quadratic forms $Q_t$ and $Q_h$ on $FG$ by setting, for $u=\sum_{g\in G}u_gg\in FG$
\begin{equation}\label{E:Q_t}
\begin{array}{lcl}
Q_t(u)&=&\sum\limits_{\{g,tg\}\subseteq G}u_gu_{tg},\\
Q_h(u)&=&\sum\limits_{\hphantom{\{}g\in G\hphantom{,tg\}}}u_gu_{hg}.
\end{array}
\end{equation}
Then each $G$-invariant quadratic form on $FG$ is a linear combination of $Q_t$'s and $Q_h$'s.

\section{Real $2$-blocks of defect zero}\label{S:defect0}

Assume that $G$ has even order, and that $B$ is a real $2$-block of $G$ which has a trivial defect group. Equivalently $B$ is a simple $F$-algebra which is a $\sigma$-invariant $FG\times G$-direct summand of $FG$. Moreover, $B$ has a unique irreducible $K$-character $\chi$ and a unique simple module $S$. %According to \cite{G} $B$ has an extended defect group of order $2$. 

Let $e_B$ be the identity element (or block idempotent) of $B$. Then
$$
e_B=e_1+e_2+\dots+e_d,
$$
where $d={\rm dim}_F(S)$ and the $e_i$ are pairwise orthogonal primitive idempotents in $FG$. Each $e_iFG$ is isomorphic to $S$. In particular $S$ is a projective $FG$-module.

Let $M$ be an $RG$-lattice whose character is $\chi$. Then $M/J(R)M\cong S$, as $FG$-modules. Now $M$ has a quadratic geometry because $\chi$ has Frobenius-Schur indicator $+1$. Thus $S$ has a quadratic geometry.

By \cite{GW93} there exists an involution $t$ in $G$ such that the restriction of the form $Q_t$ of \eqref{E:Q_t} to $e_1FG$ is non-degenerate. It follows that $e_1$ can be chosen so that $e_1=e_1^{t\sigma}=te_1^\sigma t$. We note that it can be shown that $\langle t\rangle$ is an extended defect group of $B$ and $S$ is a direct summand of $F_{C_G(t)}{\uparrow^G}$.

As $e_B=e_B^{t\sigma}$, we have $e_B=e_1+e_2^{t\sigma}+\dots+e_d^{t\sigma}$, and each $e_i^{t\sigma}$ is primitive in $FG$ and $e_1e_i^{t\sigma}=0=e_i^{t\sigma}e_1$, for $i>1$.

Suppose next that $V$ is a $B$-module, equipped with a (possibly degenerate) $G$-invariant symmetric bilinear form $\langle\,,\,\rangle$. The $G$-invariance is equivalent to $\langle ux,v\rangle=\langle u,vx^\sigma\rangle$, for all $u,v\in V$ and $x\in FG$. Now $e_1e_i=0$, for $i>1$. So
$$
\langle Ve_1,Ve_i^\sigma\rangle=0,\quad\mbox{for $i>1$.}
$$

Following \cite{KOW}, we define a bilinear form $b$ on the $F$-space $Ve_1$ by
$$
b(ue_1,ve_1):=\langle ue_1,ve_1t\rangle,\quad\mbox{for all $ue_1,ve_1\in Ve_1$.}
$$
Then $b$ is symmetric, as
$$
b(ue_1,ve_1)=\langle ue_1t,ve_1\rangle=\langle ve_1,ue_1t\rangle=b(ve_1,ue_1).
$$

Now consider the radicals of the forms
$$
\begin{array}{ll}
{\rm rad}(V)		&:=\{u\in V\mid\langle u,v\rangle=0,\forall v\in V\},\\
{\rm rad}(Ve_1)		&:=\{ue_1\in Ve_1\mid b(ue_1,ve_1)=0,\forall ve_1\in Ve_1\}.
\end{array}
$$
We include a proof of Lemma 4.5 of \cite{KOW} for the benefit of the reader:

\begin{Lemma}\label{L:KOW}
${\rm rad}(Ve_1)={\rm rad}(V)e_1$ and $Ve_1/{\rm rad}(Ve_1)\cong(V/{\rm rad}(V))e_1$.
\begin{proof}
Let $u\in{\rm rad}(V)$ and $ve_1\in Ve_1$. Then
$$
b(ue_1,ve_1)=\langle ue_1,ve_1t\rangle=\langle u,ve_1te_1^\sigma\rangle=0.
$$
So ${\rm rad}(Ve_1)\supseteq{\rm rad}(V)e_1$. Now let $ue_1\in{\rm rad}(Ve_1)$ and $v\in V$. Writing $v=\sum_{i=1}^dve_i^\sigma$, we have
$$
\langle ue_1,v\rangle=\sum_{i=1}^d\langle ue_1,ve_i^\sigma\rangle=\langle ue_1,ve_1^\sigma\rangle=b(ue_1,vte_1)=0.
$$
So ${\rm rad}(Ve_1)\subseteq{\rm rad}(V)e_1$. The stated equality follows.

We have an $F$-vector space map $\phi:Ve_1\rightarrow(V/{\rm rad}(V))e_1$ such that $\phi(ve_1)=ve_1+{\rm rad}(V)$. Now $(v+{\rm rad}(V))e_1=ve_1+{\rm rad}(V)$ as ${\rm rad}(V)e_1\subseteq {\rm rad}(V)$. So $\phi$ is onto. Moreover, ker$(\phi)={\rm rad}(V)e_1$. The stated isomorphism follows from this.
\end{proof}
\end{Lemma}

\section{Brauer characters of symmetric groups}

Let $n$ be a positive integer. Corresponding to each partition $\lambda$ of $n$, there is a Young subgroup $\Sigma_\lambda$ of $\Sigma_n$ and a permutation $R\Sigma_n$-module $M^\lambda:={\rm Ind}_{\Sigma_\lambda}^{\Sigma_n}(R_{\Sigma_\lambda})$. This module has a $\Sigma_n$-invariant symmetric bilinear form with respect to which the permutation basis is orthonormal. The {\em Specht lattice} $S^\lambda$ is a uniquely determined $R$-free $R\Sigma_n$-submodule of $M^\lambda$ c.f. \cite[4.3]{J}. Then $S^\lambda\otimes_RK$ is an irreducible $K\Sigma_n$-module and all irreducible $K\Sigma_n$-modules arise in this way.

Now $S^\lambda$ is usually not a self-dual $R\Sigma_n$-module; the dual module $S_\lambda:=S^{\lambda*}$ is naturally isomorphic to $S^{[1^n]}\otimes_RS_R^{\lambda^t}$ where $\lambda^t$ is the transpose partition to $\lambda$. Note that $S^{[1^n]}$ is the $1$-dimensional {\em sign module}.

Set $\ov{S^\lambda}:=S^\lambda/JS^\lambda$. Then $\ov{S^\lambda}$ is a Specht module for $F\Sigma_n$. It inherits an $\Sigma_n$-invariant symmetric bilinear form $\langle\,,\,\rangle$ from $S^\lambda$. This form is nonzero if and only if $\lambda$ is $2$-regular (i.e. if $\lambda$ has different parts).

Suppose that $\lambda$ is $2$-regular. Then $D^\lambda:=\ov{S^\lambda}/{\rm rad}(\ov{S^\lambda})$ is a simple $F\Sigma_n$-module, and all simple $F\Sigma_n$-modules arise uniquely in this way. The $D^\lambda$ are evidently self-dual. Indeed, $\langle\,,\,\rangle$ induces a nondegenerate form on $D^\lambda$, which by Fong's Lemma is symplectic if $D^\lambda$ is non-trivial. Note that $\ov{S^{[1^n]}}$ is the trivial $F\Sigma_n$-module, as char$(F)=2$. It follows that the dual of a Specht module in characteristic $2$ is a Specht module:
$$
\ov{S_\lambda}\cong\ov{S^{\lambda^t}}.
$$

Let $B$ be a $2$-block of $\Sigma_n$. Then $B$ is determined by an integer {\em weight} $w$ such that $n-2w$ is a nonnegative triangular number $k(k+1)/2$. The partition $\delta:=[k,k-1,\dots,2,1]$ is called the {\em $2$-core} of $B$. Each defect group of $B$ is $\Sigma_n$-conjugate to a Sylow $2$-subgroup of $\Sigma_{2w}$.

Recall that the {\em $2$-core} of a partition $\lambda$ is obtained by successively stripping removable domino shapes from $\lambda$. We attach to $B$ all partitions of $n$ which have $2$-core $\delta$.

Set $m:=n-2w$ and identify $\Sigma_{2w}\times\Sigma_m$ with a Young subgroup of $\Sigma_n$. Now $\Sigma_m$ has a $2$-block $B_\delta$ of weight $0$ and $2$-core $\delta$. This block is real and has a trivial defect group. Moreover, $S_K^\delta$ is the unique irreducible $K\Sigma_m$-module in $B_\delta$ and $D^\delta=\ov{S^\delta}$ is the unique simple $B_\delta$-module. It is important to note that $D^\delta$ is a projective $F\Sigma_m$-module and every $F\Sigma_m$-module in $B_\delta$ is semi-simple.

Let $e_\delta$ be the block idempotent of $B_\delta$. Following Section \ref{S:defect0}, choose an involution $t\in\Sigma_m$ and a primitive idempotent $e_1$ in $F\Sigma_m$ such that $e_1=e_1e_\delta$ and $e_1^{t\sigma}=e_1$. Note that dim$_F(D^\delta e_1)=1$.

Let $\mu$ be a $2$-regular partition in $B$. Regard $V:=\ov{S^\mu}e_\delta$ as an $F\Sigma_{2w}\times\Sigma_m$-module by restriction. Then $Ve_1$ is an $F\Sigma_{2w}$-module, as the elements of $\Sigma_{2w}$ commute with $e_1$. Indeed
$$
V\cong Ve_1\otimes_F D^\delta\quad\mbox{as $F\Sigma_{2w}\times\Sigma_m$-modules.}
$$
Now $\ov{S^\mu}$ and hence $V$ affords a $\Sigma_{2w}\times\Sigma_m$-invariant symmetric bilinear form $\langle\,,\,\rangle$ such that $V/{\rm rad}(V)=D^\mu e_\delta$. It then follows from Lemma \ref{L:KOW} that we may use the identity $e_1^{t\sigma}=e_1$ to construct a symmetric bilinear form $b$ on $Ve_1$. Moreover, $Ve_1/{\rm rad}(Ve_1)\cong D^\mu e_1$. So the $F\Sigma_{2w}$-module $D^\mu e_1$ inherits a nondegenerate symmetric bilinear form $b$. Reviewing the construction of $b$ from $\langle\,,\,\rangle$, we see that $b$ is $\Sigma_{2w}$-invariant (as $t\in\Sigma_m$ commutes with all elements of $\Sigma_{2w}$, and $\langle\,,\,\rangle$ is $\Sigma_n$-invariant).

\begin{Lemma}\label{L:symplectic}
Suppose that $\mu\ne[k+2w,k-1,\dots,2,1]$. Then $D^\mu e_1$ affords a non-degenerate $\Sigma_{2w}$-invariant symplectic bilinear form.
\begin{proof}
In view of Lemma \ref{L:trivial} and the discussion above, it is enough to show that $D^\mu e_1$ has no trivial $F\Sigma_{2w}$-submodules. Suppose otherwise. Then $F_{\Sigma_{2w}}\otimes_FD^\delta$ is a submodule of the restriction of $D^\mu$ to $\Sigma_{2w}\times\Sigma_m$. But $D^\mu$ is a submodule of $\ov{S_\mu}$. So $D^\delta$ is a submodule of ${\rm Hom}_{F\Sigma_{2w}}(F_{\Sigma_{2w}},\ov{S_\mu})$ as $F\Sigma_m$-modules.

We have $F$-isomorphisms
$$
\begin{aligned}{}
{\rm Hom}_{F\Sigma_{2w}}(F_{\Sigma_{2w}},\ov{S_\mu})
&\cong
{\rm Hom}_{F\Sigma_n}(M^{[2w,1^m]},\ov{S_\mu}),\hspace{0.8em}\mbox{by Eckmann-Shapiro}\\
&\cong
{\rm Hom}_{F\Sigma_n}(\ov{S^\mu},M^{[2w,1^m]}),\hspace{0.8em}\mbox{as $M^{[2w,1^m]}$ is self-dual.}
\end{aligned}
$$
As $\mu$ is $2$-regular, it follows from \cite[13.13]{J} that ${\rm Hom}_{F\Sigma_n}(\ov{S^\mu},M^{[2w,1^m]})$ has a basis of semistandard homomorphisms. The argument of Theorem 4.5 of \cite{H} now applies, and shows that
$$
{\rm Hom}_{F\Sigma_{2w}}(F_{\Sigma_{2w}},\ov{S_\mu})\cong \ov{S^{\mu^t\backslash[1^{2w}]}}\quad\mbox{as $F\Sigma_m$-modules.}
$$
Here $\mu^t\backslash[1^{2w}]$ is a skew-partition of $m$; it is empty if $\mu_1<2w$ (in which case ${\rm Hom}_{F\Sigma_{2w}}(F_{\Sigma_{2w}},\ov{S_\mu})=0$). Otherwise its diagram is the set of nodes in the Young diagram of $\mu^t$ not in the top $2w$ rows of the first column. Now $\ov{S^{\mu^t\backslash[1^{2w}]}}$ has an $F\Sigma_m$-submodule isomorphic to $D^\delta$ if and only if $S_K^{\mu^t\backslash[1^{2w}]}$ has an $K\Sigma_m$-submodule isomorphic to $S_K^\delta$, as $D^\delta=\ov{S^\delta}$, and using the projectivity of $D^\delta$.

The multiplicity of $S_K^\delta$ in $S_K^{\mu^t\backslash[1^{2w}]}$ is the number of $\mu\backslash[2w]$-tableau of type $\delta^t=\delta$ which are strictly increasing along rows and nondecreasing down columns. Suppose for the sake of contradiction that such a tableau $T$ exists.

We claim that $\mu_i\leq k-i+2$ for $i=2,\dots,k$, and $\mu_i=0$ for $i>k+1$. This is true for $i=2$, as the entries in the second row of $T$ are different. Suppose that $i\geq2$ and $\mu_{i-1}\leq k-i+3$. But $\mu_i<\mu_{i-1}$, as $\mu$ is $2$-regular. So $\mu_i\leq k-i+2$, proving our claim.

On the other hand, $\mu_i\geq\delta_i=k-i+1$, for $i=1,\dots,k$, as $\mu$ has $2$-core $\delta$. It follows that $\mu\backslash\delta$ consists of the last $\mu_1-k$ nodes in the first row of $\mu$, and a subset of the nodes $(2,k),(3,k-1),\dots,(k,2),(k+1,1)$. On the other hand, $\mu$ has $2$-core $\delta$. So $\mu\backslash\delta$ is a union of domino shapes. It follows that $T$ does not exist if $\mu\ne[k+2w,k-1,\dots,2,1]$. This contradiction completes the proof of the Lemma.
\end{proof}
\end{Lemma}

Suppose that $G$ is a finite group and that $B$ is a $2$-block of $G$ with defect group $P\leq G$. Then it is known that $[G:P]_2$ divides the degree of every irreducible Brauer character in $B$. Recall that a Brauer character in $B$ has {\em height zero} if the $2$-part of its degree is $[G:P]_2$. We now prove the main result of \cite{KOW}.

\begin{Theorem}
Let $B$ be a $2$-block of $\Sigma_n$. Then $B$ contains a unique irreducible Brauer character of height $0$.
\begin{proof}
Suppose as above that $B$ has weight $w$ and $2$-core $\delta$, and let $\theta$ be a height zero irreducible Brauer character in $B$. Then $\theta$ is the Brauer character of $D^\mu$ for some $2$-regular partition $\mu$ of $n$ belonging to $B$. 

Let $P$ be a vertex of $D^\mu$. Then $P$ is a defect group of $B$. We may assume that $P$ is a Sylow $2$-subgroup of $\Sigma_{2w}$. It is easy to show that $N_{\Sigma_n}(P)=P\times\Sigma_m$, a subgroup of $\Sigma_{2w}\times\Sigma_{m}$.

Let $B_0$ denote the principal $2$-block of $\Sigma_{2w}$. Then $B_0\otimes B_\delta$ is the Brauer correspondent of $B$ with respect to $(\Sigma_n,P,\Sigma_{2w}\times\Sigma_{m})$. So the Green correspondent of $D^\mu$ with respect to $(\Sigma_n,P,\Sigma_{2w}\times\Sigma_m)$ has the form $U^\mu\otimes D^\delta$, where $U^\mu$ is an indecomposable $\Sigma_{2w}$-direct summand of $D^\mu e_1$ which belongs to $B_0$. Moreover, $U^\mu$ is the unique component of $D^\mu e_1$ that has vertex $P$.

If $\mu=[k+2w,k-1,\dots,2,1]$ it can be shown that $U^\mu$ is the trivial $F\Sigma_{2w}$-module. Suppose that $\mu\ne[k+2w,k-1,\dots,2,1]$. Lemma \ref{L:symplectic} implies that $D^\mu e_1$ has a symplectic geometry. It then follows from the first proposition in \cite{GW95} that $U^\mu$ has a symplectic geometry. In particular dim$(U^\mu)$ is even.

Now the $2$-part of dim$(U^\mu\otimes D^\delta)$ divides $2|\Sigma_m|_2=2[\Sigma_n:P]_2$. A standard result on the Green correspondence implies that the $2$-part of dim$(D^\mu)$ divides $2[\Sigma_n:P]_2$. This contradicts the assumption that $\theta$ has height zero, and completes the proof.
\end{proof}
\end{Theorem}

\section{Acknowledgement}

B. K\"ulshammer drew my attention to the preprint \cite{KOW} during a visit to Jena in April 2011. S. Kleshchev suggested I look at the restrictions of dual Specht modules, and D. Hemmer clarified the `fixed-point functors' used to prove Lemma \ref{L:symplectic}. G. Navarro gave me the example of the group of order $72$ which has a $2$-block with $3$ real irreducible Brauer characters of height $0$.

%%%%%%%%%%%%%%%%%%%%%%%%%%%%%%%%%%%%%%%%%%%%%%%%%%%%%%%%%%%%%%%%%

\end{document}